
\input amssym.def
\input amssym

\magnification = \magstep1
\hsize=15truecm
\vsize= 21.5truecm
\hoffset=0.4 truecm
\baselineskip=12pt 
\def\tit{\medskip \noindent}
\font\big=cmbx10 scaled \magstep2
\font\little=cmr8

\hyphenation{pro-jec-tive sche-mes pro-ject MURST}
 
 \def\endproof{$\bullet$\goodbreak\medskip}
\def\smallskip{\vskip 0.3truecm}
\def\medskip{\vskip 0.5truecm}
\def\bigskip{\vskip 0.7truecm}

\def\p{\Bbb P^3}
\def\pp{\Bbb P^4}
\def\ppp{\Bbb P^5}

\centerline{\big \hbox{On quadrisecant
lines of threefolds  in $\Bbb P^5$} }
\bigskip

\centerline{EMILIA MEZZETTI 
}

\bigskip\bigskip

\qquad\qquad\qquad\qquad\qquad\qquad\qquad\qquad\qquad\qquad\qquad
{\it Dedicated to Silvio Greco}
\bigskip

{\little {

 We study smooth threefolds of $\ppp$ whose quadrisecant lines don't fill up the space. We give a complete classification of
those threefolds $X$ whose only quadrisecant lines are the lines contained in $X$. Then we prove that, if $X$ admits
\lq\lq true'' quadrisecant lines, but they don't fill up $\ppp$, then  either $X$ is contained in a cubic hypersurface, or
it contains a family of dimension at least two of plane curves of degree at least four.
 }}

\tit 
\noindent {\bf  
Introduction.}
\smallskip
The classical theorem of general projection for surfaces 
says that a general projection in $\p$ of a smooth complex projective surface $S$ 
of $\ppp$ is a surface $F$ with 
{\sl ordinary singularities} i.e.  its singular locus 
is either empty or is a curve $\gamma $ such that

\noindent\item {(i)} $\gamma$ is 
either non singular or has at most
a finite number of ordinary 
triple points;

\noindent\item {(ii)} every smooth
 point of $\gamma $ is either
a  nodal point 
 or a  pinch-point of $F$;

\noindent\item {(iii)}
the general point of
$\gamma$ is a nodal point
 for $F$;

\noindent\item {(iv)} every triple 
point of $\gamma $ is an
ordinary triple point of $F$.

\noindent(see [6], [11])

Moreover $\gamma$ is empty if and only if  $S$ is already contained in a $\p$.
\medskip
Note that the projection to $\p$ can be split in two steps:
in the first step from $\ppp$ to $\pp$ $S$ acquires only double points, while triple points appear
only in the second step from $\pp$ to $\p$. 

The problem of classifying the surfaces $S$ such that $F$
does not have any triple point is equivalent to the problem of classifying the intermediate surfaces
$S'$ of
$\pp$ whose trisecant lines don't fill up $\pp$, or \lq\lq without apparent triple points'' in the old
fashioned terminology. This problem   had been tackled by Severi in [17]. His approach was based on the
description of hypersurfaces of $\pp$ containing a
$3$--dimensional family of lines: they are quadrics and hypersurfaces birationally fibered by planes.
By consequence his theorem says that a surface $S'$ without apparent triple points either is contained in a quadric
or is birationally fibered by plane curves of degree at least $3$. Recently Aure ([1])  made this result 
precise  under smoothness assumption, proving that, if  a surface $S'$ as above is not contained in a quadric,
then it is an elliptic normal scroll.

\medskip
In the study of threefolds, several analogous questions appear, not all completely answered yet. Here we
are concerned mainly with smooth threefolds of $\ppp$ and their projections to $\pp$. We
want to study their $4$--secant lines, trying in particular to describe threefolds whose $4$--secant
lines don't fill up the space. We first study threefolds $X$ whose only $4$--secant lines are the lines
contained in $X$: we give a complete description of them (Theorem 2.1). Then we consider the
threefolds with a $5$--dimensional family of $4$--secant lines (or more generally $k$--secant lines,
with $k\geq 4$): we find that these lines cannot fill up $\pp$ and that $X$ is birationally ruled by
surfaces of $\p$ of degree $k$ (Theorem 2.3). There are no examples of this situation and it seems sensible to guess that
in fact it cannot happen. 

The general situation is that of
$3$--folds whose
$4$--secant lines form a family  of dimension four, i.e. a congruence of lines. To understand 
the case of a congruence of order $0$, i.e. of lines not filling up $\ppp$, we 
imitate
the approach of Severi: we have to look at hypersurfaces
$Y$ of
$\ppp$ covered by a
$4$--dimensional family of lines. We find that a priori there are many possibilities for such
hypersurfaces.
 More precisely, if we consider a general hyperplane section $V$ of such a $Y$, this is a threefold of
$\pp$ covered by a $2$--dimensional family of lines. The threefolds like that are studied in  
[12], where the following result is proved:

\medskip

\noindent{\bf Theorem 0.1}.
{\sl Let $V\subset\pp$ be a projective, integral hypersurface 
over an algebraically closed field  of characteristic zero,
covered by lines. Let $\Sigma\subset\Bbb G(1,4)$ denote the Fano scheme   of the lines on $V$.
 Assume that $\Sigma$ is generically reduced of dimension $2$.  
 Let $\mu$ denote the number of lines of $\Sigma$ passing through a general point of $V$ and $g$ the sectional genus of
$V$, i.e. the geometric genus of a plane section of $V$.  Then $\mu\leq 6$ and one of the following happens:

\noindent\item {(i)} $\mu=1$, i.e. $V$ is birationally a scroll over a surface;

\noindent\item {(ii)}  $V$ is birationally ruled by smooth quadric surfaces over a curve ($\mu=2$);
 
\noindent\item {(iii)} $V$ is a cubic hypersurface with singular  locus of dimension at most
one; if $V$ is smooth,
then $\Sigma$ is irreducible and $\mu=6$;

\noindent\item {(iv)} $V$ has degree $d\leq 6$,  $g=1$, $2\leq \mu\leq 4$ and $V$ is a
projection in
$\pp$ of one of the following:

\qquad   -
 a complete intersection of two hyperquadrics in $\ppp$, $d=4$; 

\qquad   -  a
section of $\Bbb G(1,4)$ with a $\Bbb P^6$, $d=5$;

\qquad   - a hyperplane section of $\Bbb P^2\times\Bbb P^2$, $d=6$;

\qquad   -  $\Bbb P^1\times \Bbb P^1\times \Bbb P^1$, $d\leq 6$.}

\medskip
\noindent To apply this result to   a fourfold $Y$ generated by the $4$--secant lines of a threefold, it
is necessary first of all to understand the meaning of the assumption of generic reducedness on
$\Sigma$.  
We prove that this hypothesis is equivalent to the non--existence of a fixed tangent plane to $V$
along a general line of $\Sigma$. Threefolds $V$ not satisfying this assumption are then 
described in Proposition 3.4. 

So it is possible to perform an analysis of the possible cases for the fourfold $Y$. This leads to  a
result very similar to the theorem of Severi for surfaces quoted above:

\medskip
\noindent{\bf Theorem 0.2}.

{\sl Let $X$ be a smooth non--degenerate threefold of $\ppp$ not contained in a quadric. Let $\Sigma$ be an
irreducible component of dimension $4$ of $\Sigma_4(X)$ such that a general line of $\Sigma$ is
$k$--secant $X$ ($k\geq 4$). Assume that the union of the lines of $\Sigma$ is a hypersurface $Y$. Then
either $Y$ is a cubic or $Y$ contains a family of planes of dimension $2$ which cut on $X$ a family of
plane curves of degree $k$.}
\medskip

Recently, a different approach to the study of multisecant lines of smooth threefolds of $\ppp$ has been
considered by Sijong Kwak ([8]). It is based on the well--known monoidal construction. He proves that,
if the $4$--secant lines of $X$ don't fill  up $\ppp$, then either $h^2({\cal O}_X)\not=0$ or $h^1({\cal
O}_X(1))\not=0$. Moreover he gives an explicit formula for  $q_4(X)$, the number of $4$--secant lines
through a general point of $\ppp$, depending on $\deg X$, on the sectional genus and on the two Euler
characteristics $\chi({\cal O}_X)$ and $\chi({\cal O}_S)$, where $S$ is  a general hyperplane section.
It is interesting to note that, testing this formula on all known smooth threefolds of $\ppp$, one
gets $q_4(X)=0$ only for those contained in a cubic hypersurface.
\bigskip

{\bf Aknowledgment.} I warmly thank Silvio Greco and Dario Portelli for several dicussions 
on the topic of this paper and
for encouragement. I would like to thank Prof. Sijong Kwak for inviting me to KIAS of Seoul. There I
benefitted also from interesting conversations with Prof. Fjodor Zak, whom I gratefully thank. 

This work has been done in the framework of the activities of the GNSAGA - INdAM project ``Cohomological
properties of projective schemes" and of the MURST project \lq\lq Geometria sulle variet\`a
algebriche''. The author is member of EAGER.

\tit 
\noindent {\bf 1.  
Multisecant 
lines of threefolds in $\ppp$.}
\smallskip
Let $X$ be an integral 
smooth threefold 
of $\Bbb P^5$ not contained in a hyperplane. To define the multisecant lines of $X$, we follow the
approach of Le Barz ([9]).  Let
$k\geq 2$ be an integer number. Let
$Hilb^k{\Bbb P^5}$ be
 the Hilbert scheme 
of subschemes of length $k$
 of ${\Bbb P^5}$, and
$Hilb^k_c  {\Bbb P^5}$ be
 its smooth open 
subvariety parametrizing curvilinear
subschemes, i.e. subschemes which are contained 
in a smooth curve. Let
 $Al^k {\Bbb P^5}$
denote the subscheme of 
$Hilb^k_c {\Bbb P^5}$ of  
length $k$ subschemes lying on a line
and $Hilb^k_cX$ that 
of subschemes contained in $X$.
The following cartesian
 diagram defines $Al^k X$, the scheme of aligned $k$-tuples
of points of $X$:
$$\matrix{ Al^k X &  
\longrightarrow & Al^k\Bbb P^5
 \cr
 \downarrow &     &       
\downarrow \cr
Hilb^k_cX  & \longrightarrow 
 & Hilb^k_c\Bbb P^5.\cr} 
$$

We have: $\dim Hilb^k_c 
{\Bbb P^5}=5k$, $\dim  
Al^k{\Bbb P^5}=10 + (k-2)=8+k$, 
$\dim Hilb^k_cX=3k$;
 so, if $Al^k X$ is non--empty, then any irreducible
component  of   its
has dimension at least $(8+k)+(3k)-(5k)=8-k$.

\par Let now $$a\colon Al^k
\Bbb P^5\longrightarrow \Bbb G(1,5)$$
 be the
 natural
 map (axe) to the Grassmannian 
of lines of $\Bbb P^5$.  Note that all
fibers of
$a$  have dimension $k$.

The image scheme $\Sigma_k(X):=a(Al^k(X))$
is by definition the family of $k$-secant lines of $X$.
Clearly  all lines contained in $X$ belong to
$\Sigma_k(X)$. If $Al^k X=\emptyset$, then obviously also
$\Sigma_k(X)=\emptyset$: in this case no line cuts $X$ in
at least $k$ points or is contained in $X$.

Let us consider now the restriction $ \bar a$ of 
$a$ to an irreducible component $\Sigma$ of
 $Al^k X$:
$$ \bar a \colon \Sigma\longrightarrow 
\bar a (\Sigma)\subset \Bbb G(1,5).$$
 
We have: 
$\dim \bar a(\Sigma)=\dim \Sigma -\dim \Sigma_l$,
 where $\Sigma_l:=\bar a^{-1}(l)$ is the fibre over
 $l$, a general line of $\bar a(\Sigma)$.
There are two possibilities, i.e.
either $\dim\Sigma_l=k$ if $l\subset X$, 
or
dim $\Sigma_l=0$ if $l\cap X$ is a
 scheme of finite length. By consequence, either
dim $\bar a(\Sigma)=\dim \Sigma -k$, if
 any line of $\bar a(\Sigma)$ is contained in $X$,
 or else $\dim \bar a(\Sigma)=
\dim \Sigma$ if a general line of
 $\bar a(\Sigma)$ 
is not contained in $X$.

\medskip
Some rather precise information on the families of $k$--secant lines  of threefolds
for particular $k$ come from the classical theorems of \lq\lq\thinspace general projection''. For smooth
curves in
$\p$ and smooth surfaces in
$\pp$ there are  very precise theorems, describing the singular locus of the projected variety (see
[7], [11], [1]).

{}From these results, passing to  general sections with linear spaces of
dimension $3$ and
$4$, it follows that a general projection $X'$ in $\pp$ of a smooth threefold $X$ of $\ppp$ acquires a
double  surface $D$, i.e. a surface whose  points have multiplicity at least two on $X'$, and a 
triple curve $T\subset D$, i.e. a curve whose  points have multiplicity at least three on $X'$.
Moreover, $D$ is non--empty
unless
$X$ is degenerate and $T$ is non--empty unless
$X$ is contained in a quadric. In terms of multisecant lines, this means that, through a general point
$P$ of
$\ppp$, there passes a
$2$--dimensional family of   $2$- secant lines of $X$: we have that $\Sigma_2(X)$ is irreducible of
dimension
$6$. If moreover $X$ is not contained in a quadric, then the trisecant lines through $P$ form a family
of dimension $1$, so $\Sigma_3(X)$ has dimension $5$ and its lines fill up $\ppp$. I would like to
emphasize that $D$ is truly double and $T$ is truly triple for $X'$, or, in other words, a general 
secant line of $X$ is not trisecant and a general trisecant line  is not quadrisecant.

On the other hand, it has been proved that $X'$ does not have any point of multiplicity $5$ or more
(see [15], [10]). Hence the $5$--secant lines of $X$ never fill up $\ppp$.
\bigskip
{\bf Remark 1.1.} It is interesting to note that no smooth threefold $X$ in $\ppp$ has $Al^3
(X)=\emptyset$. Indeed, if so, a general curve section $C$ of $X$ would be a smooth curve of $\p$
without trisecant lines. It is well known that such a curve $C$ is either a skew cubic or an elliptic
quartic. So $X$ could be either $\Bbb P^1 \times \Bbb P^2$ or a complete intersection of two quadrics: 
in both cases, $X$ is an intersection of quadrics, so the trisecant lines are necessarily
 contained in $X$. But both threefolds contain lines: they form a family of
dimension
$3$ in the first case and of dimension
$2$ in the second one.

\tit 
\noindent {\bf 2.  Quadrisecant lines: special cases. 
}
\smallskip

The first case we consider is that of threefolds without \lq\lq\thinspace true'' quadrisecant lines.
\smallskip

\noindent{\bf Theorem 2.1}

{\sl Let $X$ be a smooth threefold of $\ppp$. Then $Al^4(X)\not =\emptyset$. If  its quadrisecant lines
are all contained in $X$, then $\sigma_4:=\dim\Sigma_4(X)\leq 4$ and one of the following possibilities
occurs:

$\sigma_4=4$, $X$ is a $\p$;

$\sigma_4=3$, $X$ is a quadric hypersurface (contained in a hyperplane of $\ppp$), or $\Bbb P^1 \times
\Bbb P^2$;

$\sigma_4=2$, $X$ is a cubic hypersurface (contained in a hyperplane of $\ppp$), or a complete
intersection of type
$(2,2)$, or a Castelnuovo threefold, or a Bordiga scroll;

$\sigma_4=1$, $X$ is a complete intersection of type $(2,3)$, or an inner projection of a complete
intersection of type $(2,2,2)$ in $\Bbb P^6$;

$\sigma_4=0$, $X$ is a complete intersection of type $(3,3)$.}
\smallskip
\noindent{\bf Proof.} 
By [16],  the maximal dimension of a family of lines contained in a threefold $X$ is
$4$, and the maximum is attained only by linear spaces. Moreover, if the dimension is $3$, then
either $X$ is  a quadric or it is birationally a scroll over a curve. Being $X$ smooth, in the last
 case  $X$ is $\Bbb P^1 \times \Bbb P^2$ (see [14]).

If $\sigma_4\leq2$, then a general hyperplane section $S$ of $X$ contains only a finite number of lines
and does not possess any other $4$--secant line. In [3] one proves that there is a finite explicit
list of such surfaces $S$. They have all degree at most $9$ and are all arithmetically Cohen-Macaulay, except
the elliptic scroll.
The smooth threefolds $X$ having them as general hyperplane sections are all described (see for
instance [5]) and are
precisely those appearing in the list above. More precisely, a Castelnuovo threefold has degree $5$, 
its ideal is generated
by the maximal minors of a $2\times 3$ matrix of forms: in the first two columns the entries are linear
while in the third one they are quadratic, $X$ is fibered by quadrics over $\Bbb P^1$. The Bordiga
scroll has degree $6$, its ideal is generated by the maximal minors of a $3\times 4$ matrix of linear
forms, it is a scroll over $\Bbb P^2$.
Finally, the computation of the dimension of the family of lines contained in a smooth complete
intersection as above is classical.  Note that in all cases $X$ contains lines, so $Al^4(X)\not
=\emptyset$.
\endproof

{\bf Remark 2.2.}
Note that all threefolds whose only quadrisecant lines are the lines contained in them are cut out by
quadrics and cubics.

\bigskip
{}From now on we will consider only smooth non--degenerate threefolds in $\ppp$ such that the general
line of at least one irreducible component of $\Sigma_4(X)$ is not contained in $X$. Hence the
dimension of such a component $\Sigma$ is at least $4$. On the other hand $\dim\Sigma<6$, otherwise
every secant line would be quadrisecant, which is excluded by general projection theorems.
If the dimension of such a component is $5$, then we have the following result.
\medskip
\noindent {\bf Theorem 2.3}.

{\sl Let $X$ be a smooth non--degenerate threefold of $\ppp$, let $\Sigma$ be an irreducible component of
$\Sigma_4(X)$ of dimension $5$. Then the lines of $\Sigma$ don't fill up $\ppp$. More precisely
either their union is 
 a quadric or it is a hypersurface birationally ruled by $\p$'s over a curve. }
\smallskip
\noindent{\bf Proof.}
Let $H$ be a general hyperplane and let $S:=X\cap H$, $\Sigma':=\Sigma\cap \Bbb G(1,H)$. $S$ is a
smooth surface of $\pp$ and $\Sigma'$ is a family of dimension $3$ of quadrisecant lines of $S$. From
the general projection result for surfaces, it follows that the lines of $\Sigma'$ don't fill up $H$, so
their union is a hypersurface $V$ in $H$. By [16], either $V$ is a quadric or it is birationally
fibered by planes. In the first case, $V$ lifts to a quadric containg $X$ and all its quadrisecant
lines. 

In the second case, the planes of $V$ cut on $S$ a one--dimensional family of plane curves of
degree, say, $a$: since the lines of these planes have to be $4$--secant $S$, then $a\geq 4$. Coming
back to $\ppp$, $X$ contains a family of dimension at least $4$ of plane curves of degree at least
$4$. Let
$W$ be the subvariety of $\Bbb G(2,5)$ parametrizing those planes. We consider the focal locus of the
family $W$ on a fixed plane $\pi$ (see [4] for generalities about the theory of foci): it must contain
the plane curve of $X$ lying on $\pi$. But the matrix representing the characteristic map of $W$
restricted to
$\pi$ is a $3\times 4$ matrix of linear forms on $\pi$, so it cannot degenerate along a curve of degree
strictly bigger than $3$, unless it degenerates everywhere on $\pi$. So all planes of the family are
focal planes. Let $f$ be the projection from the incidence correspondence of $W$ to $\ppp$: the
differential of $f$ has  always a kernel of dimension two
and image of dimension
$4$. By the analogous of Sard's theorem, it follows that the union of the planes of $W$ is a variety $Y$
of dimension $4$. By
[16], we conclude that 
$Y$ is birationally ruled by $\p$'s over a curve.
\endproof

\smallskip
{\bf Remark 2.4.} Under the assumption of  Theorem 2.3, if $X$ is not contained in a quadric, then it is covered by a
one--dimensional family of surfaces of
$\p$'s of degree at least $4$, whose hyperplane sections are the plane curves covering $S$. So the plane
curves on
$X$ are cut by the planes of the $\p$'s of $Y$.

\tit 
\noindent {\bf 3.  Quadrisecant lines not filling up the space. 
}
\smallskip
We assume now that $X$ is a non--degenerate smooth threefold in $\ppp$, such that
all irreducible components of $\Sigma_4(X)$, corresponding to lines not all contained in $X$, have
dimension
$4$. A subscheme of dimension $4$ of $\Bbb G(1,5)$ is called a congruence of lines. To a congruence of
lines $\Sigma$ one associates an integer number, its {\sl order}: the number of  lines of
$\Sigma$ passing through a general point of
$\ppp$. More formally, it is the intersection number of $\Sigma$ with the Schubert cycle of lines
through a point. The order of
$\Sigma_4(X)$ will be denoted by
$q_4(X)$. It is clear that if
$X$ is contained in a quadric or in a cubic hypersurface, then this hypersurface contains also the
quadrisecant lines of
$X$, hence 
$q_4(X)=0$. It is natural to try to reverse this implication, so one can consider the following 
\smallskip
{\bf Question}. {\sl Do
there exist smooth threefolds $X$ in
$\ppp$, not contained in a cubic, but such that the $4$--secant 
lines of $X$ form a congruence with $q_4(X)=0$?}

\smallskip
{}From now on, we assume that   $H^0({\cal I}_X(3))=(0)$, $\dim
\Sigma_4(X)=4$ and $q_4(X)=0$. Let $Y$ be the hypersurface of $\ppp$ union of the $4$--secant lines of
$X$. Let $\Sigma$ be the Fano scheme
of lines contained in $Y$:  $\Sigma_4(X)$ is a union of one or more irreducible components of
$\Sigma$. Let now $H$ be a general hyperplane, $S:=X\cap H$ and $V:=Y\cap H$. So $\Sigma':=\Sigma\cap
\Bbb G(1,H)$ is the Fano scheme of lines contained in $V$ and $\Sigma_4(S)=\Sigma_4(X)\cap\Bbb G(1,H)$
parametrizes $4$--secant lines of $S$.
\smallskip

In order to apply Theorem 0.1 to our situation, we want to give  some characterization of threefolds
covered by lines with non--reduced associated Fano scheme. First of all we recall a result from [12].

Let $V$ be a threefold of $\pp$ covered by a two dimensional 
family of lines and let $\bar \Sigma$ be an irreducible component of dimension two of its Fano scheme of
lines. Let $r$ be a line on $V$ which is a general point of
$\bar\Sigma$, let $P$ be a general point of $r$ and let $\Bbb P(T_PV)$ be the projective plane
obtained by projectivization from the tangent space to
$V$ at $P$, its points correspond to tangent lines to $V$ at $P$. Choose homogeneous coordinates in
$\pp$ such that
$P=[1,0,\ldots,0]$ and
$T_PV$  has equation
$x_4=0$. In the affine chart $x_0\not=0$ with non--homogeneous coordinates $y_i=x_i/x_0$,
$i=1,\ldots,4$, $V$ has an equation
$G=G_1+G_2+G_3+\dots+G_d=0$, where the $G_i$ are the homogeneous components of $G$ and $G_1=x_4$. It is
convenient to write
$G_i=F_i+y_4H_i$, where the $F_i$ are polynomials in $y_1, y_2, y_3$. The equations $y_4=F_2=0$ (resp. 
$y_4=F_2=F_3=0$) represent lines in 
$\Bbb P(T_PV)$ which are at least $3$--tangent (resp. $4$--tangent) to $V$ at $P$.

\smallskip
\noindent{\bf Proposition 3.1.}

{\sl With the notations just introduced, $\bar\Sigma$ is reduced at  $r$ if and only
if in $\Bbb P(T_PV)$ the intersection of the conic $F_2=0$ with the cubic $F_3=0$ is reduced at the
point corresponding to $r$.}
\smallskip
\noindent{\bf Proof}.
[12], Proposition 1.3.
\endproof

In the following characterization, we need  again the notion of focal scheme of a family of lines (see
 [4]).

\smallskip

\noindent{\bf Proposition 3.2.}

{\sl Let $V$ be a threefold of $\pp$ covered by a two dimensional 
family of lines. Let $\bar \Sigma$ be an irreducible component of dimension two of the Fano scheme of
lines on $V$. Then the following are equivalent:

\noindent\item{(1)} $\bar \Sigma$ is non-reduced;

\noindent\item{(2)} $V$ has a fixed tangent space of dimension at least two along a general line of
$\bar \Sigma$;

\noindent\item{(3)} on each general line of the family $\bar \Sigma$ there is at least one focal point.}

\smallskip

\noindent{\bf Proof.}
$(1)\Leftrightarrow (2)$ 
One implication is Proposition 1.5 of [12]. This implication and the inverse one, which
is similar, follow from a local computation and from Proposition 3.1. 

$(2)\Leftrightarrow (3)$
 Let the line $r$ be a smooth, general
point of $\bar\Sigma$ and let $[x_0,\ldots ,x_4]$ be
homogeneous coordinates in
$\pp$ such that  $r$ has equations
$x_2=x_3=x_4=0.$  We consider the restriction to $r$ of the global
characteristic map relative to the family of lines $\bar\Sigma$:
$$
\chi(r): T_r\bar\Sigma\otimes{\cal O}_r\to{\cal N}_{r/\pp}.
$$
Since $T_r\bar\Sigma\otimes{\cal O}_r\simeq{\cal O}_r^2$ and 
${\cal N}_{r/\pp}\simeq{\cal O}_r(1)^3$,  the map $\chi(r)$  
can be represented by a suitable $3\times 2$ matrix $\cal M$, with linear 
entries $l_{ij}(x_0,x_1).$
If there is a fixed tangent plane $M_r$ to $V$ along $r$, it gives a (fixed) normal direction to 
$r$ in $\pp.$ If $\Lambda\subset K^5$ is the vector space of dimension two
corresponding to $r$, this normal direction can be represented by a 
vector $v\in K^5/\Lambda ,$ with $v\neq 0.$ Moreover, for any $P\in r,$ the
columns of $\cal M$ evaluated at $P$ are elements of $K^5/\Lambda.$

With  this set-up we can rephrase the condition that the tangent spaces 
to $V$ at the points of $r$ all contain the plane $M_r$ as follows, where
$v=(v_1, v_2, v_3)$:

$$
det\left(\matrix{
v_1&l_{11}(P)&l_{12}(P)\cr
v_2&l_{21}(P)&l_{22}(P)\cr
v_3&l_{31}(P)&l_{32}(P)}\right) =0
\leqno (*)$$

\smallskip\noindent for every $P\in r.$ 
The development of the above determinant is a quadratic form in $x_0,x_1$,
whose three coefficients linearly depend on $v_1,v_2,v_3$.
Since the determinant vanishes for each choice of $x_0,x_1$, these coefficients have to be identically
zero. This can be interpreted as a homogeneous linear system of three equations which admits the
non-trivial solution $(v_1,v_2,v_3)$. The determinant of the matrix of the coefficients of the system is
therefore zero. It
is a polynomial $G$,  homogeneous of degree $6$ in the coefficients of
the linear forms $l_{ij}$, which can be explicitly written. If $\varphi_{12}$, $\varphi_{13}$,
$\varphi_{23}$ are the quadratic forms given by the $2\times 2$ minors of $\cal M$,
it is possible to verify that  the resultant of any two of them is a multiple of $G$. Being $G=0$, it
follows that the polynomials $\varphi_{ij}$'s have a common linear factor.
Hence on a general $r\in\Sigma_1$  there exists a focal point.

The inverse implication is similar: if the polynomials $\varphi_{ij}$ have a common linear factor $L$,
such that
$\varphi_{ij}=L\psi_{ij}$, for all $i,j$, then the (*) takes the form $v_1\psi_{23}- v_2\psi_{13}+v_3\psi_{12}=0$: this is
an equation in $v_1, v_2, v_3$ which certainly admits a non-zero solution. This gives a vector $v\in K^5/\Lambda$, 
hence a normal direction to $r$ that generates the required plane $M_r$.
\endproof

\medskip

\noindent{\bf Proposition 3.3.}

{\sl If the equivalent conditions of Proposition 3.2 are satisfied, let $F$
 be the focal scheme on $V$. Then $F$ is a
point or a curve  or a surface. In the first case $V$ is a cone, in the second case $F$ is a
fundamental curve for the lines of $\bar\Sigma$ and $V$ is a union of cones with vertex on $F$, in the
third case
all lines of $\bar\Sigma$ are tangent to $F$.}
\smallskip
\noindent{\bf Proof.}
Let $I\subset\bar\Sigma\times\pp$ be the incidence correspondence,
and let $f: I\to V$ and  $q:I\to\bar\Sigma$ be the projections.  The focal scheme on $V$ can be seen
as the branch locus of the map $f$, i.e. the image of the ramification locus $\cal F$, which is a
surface. So $\dim F\leq 2$. 

The first two cases are clear. We have  to show that, if $F$ is a surface,
then all lines of
$\bar\Sigma$ are tangent to $V$. Let $P$ be a focal point on $r$
and assume that $P$ is a smooth point for $F$. 
 Let
$s\subset I$ be the fibre of $q$ over the point representing
$r$. The tangent space to $I$ at $(P,r)$ contains 
the tangent space to $\cal F$ at $(P,r)$, the line $s$ and the kernel of the
differential map ${\rm d}f$ of $f$ at $(P,r)$. Since $F$ is smooth at $P$, this 
latter space is transversal to $T_{(P,r)}{\cal F}$, and the image of ${\rm d}f$
is ${\rm d}f(T_{(P,r)}{\cal F})=T_PF.$ But also $s$ is transversal to
$ker({\rm d}f)$, hence $r={\rm d}f(s)\subset T_PF.$
\endproof

{\bf Remark 3.4.} 

1. One can prove that, if on each line $r$ of $\bar \Sigma$ there is also a second
focal point, possibly coinciding with the first one, then the tangent space to $V$ is fixed along $r$
and
$\bar\Sigma$ is the family of the fibres of the Gauss map of $V$ (see [13]).
In this case, clearly, only one line of $\bar\Sigma$ passes through a general point of $V$.

2. Also in the last case of Proposition 3.3, i.e. if the focal locus on $V$ is a surface $F$ and on a
general line
$r$ of
$\bar\Sigma$ there is only one simple focus, we can conclude that
only one line of $\bar\Sigma$ passes through a general point of $V$.
Indeed, first of all let us exclude that there are two lines $r$ and $r'$ of $\bar\Sigma$ which are
both tangent to $F$ at a general point $P$. Otherwise $r$ and $r'$ are both  contained in $T_PF$ and
the hyperplanes which are tangent to $V$ along $r$ vary in the pencil containing the fixed plane
$M_r$, which coincides with $T_PF$ in this case. So the pencil would be the same for $r$ and $r'$, and
 every hyperplane in the pencil would be tangent to $V$ at two points, one on $r$ and the other on $r'$,
which is impossible. So only one line of $\bar\Sigma$ passes through a general focal point on $V$. But
then a fortiori the same conclusion holds true also for a general non--focal point of $V$.

\smallskip
\noindent We are now able to prove Theorem 0.2 stated in the Introduction.

\medskip
\noindent{\bf Proof} of Theorem 0.2.

\noindent Let $V=Y\cap H$, where $H$ is a general hyperplane. Hence $V$ is a hypersurface of $\pp$
covered by a $2$--dimensional family of lines: this is the situation of Theorem 0.1. If one irreducible
component
$\bar\Sigma$ of the Fano scheme of lines on $V$ is non--reduced, then it follows from Proposition 3.3
and the subsequent Remark 3.4 that
$V$ is a cone, or a union of cones with vertices on a curve $C$, or a union of lines all tangent to a
surface $F$: in this last case only one line of $\bar\Sigma$ passes through a general point of $V$.
It is easy to check that,  in the first two cases, to have such a
$V$ as general hyperplane section,
$Y$ has to be a cone over $V$. In the third case, the lines through a general point of $Y$  form a surface which
intersects the general hyperplane $H$ in one line (Remark 3.4), so this surface is necessarily  a plane. In any
event
$Y$ contains a
$2$--dimensional family of planes, cutting plane curves on
$X$.

Now we assume  that all irreducible components
$\bar\Sigma$ of the Fano scheme of lines on $V$ are reduced. If $V$ is as in case $(i)$ of Theorem 0.1, 
i.e. if $\mu=1$, then the lines of $Y$ through a general point 
form a plane, and we are done. 

We consider now case $(ii)$:
we prove first that
$Y$ cannot be birationally fibered by smooth quadric surfaces. Assume, by contradiction, that
$Y$  contains such a family of quadrics and let $P$ be a fixed general point of $Y$.
Then only one quadric $F_P$ of the family passes through $P$, so  the lines
contained in
$Y$ and passing through  $P$  form a quadric cone $Q_P$, the intersection of $F_P$
with its tangent space at $P$. The linear span
$\p_P:=<Q_P>$ is the tangent space 
to $F_P$ at $P$. We consider the curve $C_P:=X\cap Q_P$: it is a $k$--secant
curve on the cone $Q_P$, so $\deg C_P=2k$ and $p_a(C_P)=(k-1)^2$. On the other hand $X\cap \p_P$ is 
a connected curve of degree $d=\deg X$. If it contains also another curve $C'_P$ different from $C_P$,
then every point of
$C_P\cap C'_P$ is singular for $X\cap \p_P$, so, being $X$ smooth, $\p_P$ has to be tangent to $X$ at
each point of $C_P\cap
C'_P$. But $\{\p_P\}_{P\in Y}$ is a family of dimension $4$ of $3$--spaces and the tangent spaces to
$X$ form a family of dimension $3$. Therefore every $\p_P$ should be tangent to infinitely many
quadrics of $Y$, i.e. to all quadrics of $Y$, which is impossible. So $X\cap \p_P=C_P$, $d=2k$ and the
sectional genus of $X$ is $(k-1)^2=({d\over 2}-1)^2$. But this is the Castelnuovo bound, so every curve
section of $X$ with a $3$--space is contained in a quadric, which implies that also $X$ is contained in
a quadric hypersurface: this gives the required contradiction.
As a consequence, if $V$ is as in $(ii)$ of Theorem 0.1, then $Y$ is birationally fibered by quadrics
of rank at most $3$. So the $k$--secant lines of $X$ are necessarily cut by the planes contained in
these quadrics.

It remains to analyze the four cases of $(iv)$ in Theorem 0.1, with $g=1$. If $V$ is a projection of a
complete intersection of type $(2,2)$, then also 
$Y$ is
a projection of a fourfold $Z$ of degree $4$ in $\Bbb P^6$, complete intersection of two quadrics. We
have the following diagram:
$$
\matrix{
{}&{}&\ Z&\subset&\Bbb P^6\cr
{}&{}&\pi\downarrow&&\cr
X&\hookrightarrow&\ Y&\subset& \ppp}
$$
\noindent where $\pi$ is the projection from a suitable point $P$.  $P\notin Z$, because $d=4$, hence the singular
locus of $Y$ is a threefold $D$ of degree $2$, according to the formula $\deg D=(d-1)(d-2)/2-g$, where $d=\deg Z$ and $g$
is the sectional genus, so D does not contain
$X$. Therefore the restriction of $\pi: \pi^{-1}(X)\to X$ is regular and birational: but  $X$,
being smooth, is linearly normal, so  $\pi^{-1}(X)$ is already contained in a $\ppp$ and the
projection is an isomorphism. In this  case $\deg X<\deg Z=4$, but the smooth threefolds of low degree in $\ppp$ are
all completely described (see for instance [2]) and this possibility is excluded. 

The second possibility for $V$ is being a projection of $\Bbb G(1,4)\cap\Bbb P^6$  of degree $5$. So
$Y$ is
a projection from a line $\Lambda$ of a fourfold $Z$ of degree $5$ in $\Bbb P^7$. Arguing as in the
previous case, we get that  $\Lambda\cap Z=\emptyset$, then either $X$ is contained in the double
locus of $Y$, which has degree $5$, or $\pi^{-1}(X)$ is contained in a $\ppp$ and again  $\deg X<5$.
Both possibilities are excluded as before. 

The last case is when $V$ is a projection of a threefold of degree $6$ and sectional genus one of $\Bbb
P^7$. If $\Lambda\cap Z=\emptyset$, it can be treated in the same way, observing that in this case the degree of the double
locus of
$Y$ is $9$. So $\Lambda\cap Z\not=\emptyset$ and the intersection should contain the whole centre of projection. But then
$\deg Y=3$.
\endproof

\bigskip

{\bf References}

\noindent\item {[1]}  Aure,  A.B.: {\it
 {The smooth surfaces in
$\pp$ without apparent triple points}}, Duke Math. J.
 {\bf 57} (1988), 423-430

\vskip 0.1 cm

\noindent\item {[2]} Beltrametti, M., Schneider, M., Sommese, A.J.
{\it Threefolds of degree 9 and 10 in ${\Bbb P}^5$}, Math. Ann. {\bf 288}, No.3 (1990), 413-444

\vskip 0.1 cm

\noindent\item {[3]} Bertolini, M., Turrini, C.: {\it Surfaces in ${\Bbb P}^4$ with no
quadrisecant lines}, Beitr\" age Algebra Geom. {\bf 39}, no. 1 (1998), 31--36. 

\vskip 0.1 cm
\noindent\item {[4]}  Chiantini, L.,  Ciliberto, C.: {\it {A few remarks on the lifting problem}}, Journ\'ees de
G\'eom\'etrie Alg\'ebrique d'Orsay (Orsay, 1992), Ast\' erisque No. 218 (1993), 95--109

\vskip 0.1 cm  

\noindent\item {[5]} Decker, W., Popescu, S.: {\it {On surfaces in ${\Bbb P}^4$ and 3-folds in
${\Bbb P}^5$}},  Vector bundles in algebraic geometry. Proceedings of the 1993 Durham symposium, 
Cambridge University Press., Lond. Math. Soc. Lect.
Note Ser. 208, 69-100 (1995)
\vskip 0.1 cm

\noindent\item {[6]}   Griffiths, P., Harris, J.: 
\lq\lq\thinspace
Principles of Algebraic Geometry", J. Wiley \& Sons, New York - 
Chichester - Brisbane - Toronto, 1978

\vskip 0.1 cm

\noindent\item {[7]}   Hartshorne, R.:
\lq\lq\thinspace
Algebraic Geometry", Springer, Berlin--Heidelberg--New York, 1977

\vskip 0.1 cm 

\noindent\item{[8]} Kwak, S.: {\it Smooth threefolds in $\ppp$ without apparent triple or quadruple
points and a quadruple--point formula}, Math. Ann. {\bf 320}, no. 4 (2001), 649--664  

\vskip 0.1 cm

\noindent\item {[9]} Le Barz, P.: {\it Platitude et non--platitude de certains sous--schemas de
Hilb$^k\Bbb P^n$}, J. Reine Angew. Math. {\bf 348} (1984),  116-134

\vskip 0.1 cm
\noindent\item {[10]} Mather, John N.: {\it Generic projections},  Ann. of Math. (2) {\bf 98} (1973),
226--245

 \vskip 0.1 cm

\noindent\item {[11]}  Mezzetti, E.,  Portelli, D.:
{\it {A tour through some classical theorems on algebraic surfaces}}, An. \c Stin\c t. Univ. Ovidius
Constan\c ta Ser. Mat. {\bf 5} no. 2, (1997),  51--78

\vskip 0.1 cm 

\noindent\item {[12]}  Mezzetti, E.,  Portelli, D.:
{\it {On threefolds covered by lines}}, Abh. Math. Sem. Univ. Hamburg {\bf 70} (2000), 211 - 238

\vskip 0.1 cm

\noindent\item {[13]} Mezzetti, E., Tommasi, O.: {\it On projective varieties of dimension $n+k$ covered
by $k$-spaces}, preprint  (2001), to appear in Illinois J. Math.

\vskip 0.1 cm

\noindent\item {[14]} Ottaviani, G.: {\it On $3$-folds in $\ppp$ which are scrolls}, Ann. Scuola
Norm. Sup. Pisa Cl. Sci. (4) {\bf 19} no. 3,
 (1992),  451--471

\vskip 0.1 cm

\noindent\item {[15]}  Ran, Z.:  {\it {The
(dimension $+2)$-secant Lemma}}, Invent. math.
 {\bf 106} (1991), 65-71

\vskip 0.1 cm

\noindent\item {[16]}  Segre, B.:
 {\it {Sulle $V_n$ contenenti pi\`u di
$\infty^{n-k}$ $S_k$}}, I, Lincei - Rend. Sc. fis. mat. nat. {\bf 5} (1948), 193 - 197, II, 
 Lincei - Rend. Sc. fis. mat. nat. {\bf 5} (1948), 275 - 280 

\vskip 0.1 cm

\noindent\item {[17]}   Severi,  F.:
 {\it {Intorno ai punti doppi
improprii di una superficie 
generale dello spazio a quattro 
dimensioni, e ai suoi punti tripli 
apparenti}}, Rend. Circolo 
Matematico di Palermo {\bf 15}
(1901), 33-51

\bigskip
\baselineskip=0.4truecm

{\little
{\it\noindent Dipartimento 
di Scienze Matematiche

\noindent Universit\`a di Trieste

\noindent  Via Valerio 12/1 
 
\noindent 34100 Trieste, ITALY

\noindent e-mail: }
mezzette@univ.trieste.it
}

\end